\documentclass[12pt]{ article}
\usepackage{amstext,amsfonts,amsmath,amsthm,graphicx,amssymb,amscd,epsfig}

\newtheorem{Theorem}{Theorem}

\newtheorem{Proposition}{Proposition}
\newtheorem{Lemma}{Lemma}

\begin{document}

\def\BC{\hbox{\rm\ l\hskip -2.2truemm C}}
\def\C{\hbox{\rm\ l\hskip -1.5truemm C}}
\def\P{\hbox{\rm\ l\hskip -0.7truemm P}}
\def\N{\hbox{\rm\bf N}}
\def\Z{\hbox{\rm\bf Z}}
\def\mult{{\rm mult}}
\def\ord{{\rm ord}}

\author{ Nguyen Van Chau\thanks {  Supported in part by the National Basic Program on Natural Science, Vietnam.} \\ {\small Institute of Mathematics, 18 Hoang Quoc Viet Road, 10307 Hanoi,  Vietnam.} \\  {\small E-mail:{\it nvchau@math.ac.vn}} }

\title{Integer points on a curve and  the plane Jacobian problem\footnote {Submitted to J. Ann. Polon. Math. , June, 2005}}  



\date {}

\maketitle

\begin{abstract}
 A polynomial map $F=(P,Q)\in \Z [x,y]^2$ with Jacobian $JF:=P_xQ_y-P_yQ_x\equiv 1$ has polynomial inverse of integer coefficients if the complex plane curve $P=0$ has infinitely many integer points.

2000{ \it Mathematical Subject Classification:} 14R15.
  
\end{abstract}

\medskip

\noindent {\bf 1. Introduction.} Let $F:\C^n\longrightarrow \C^n$ be a polynomial map with integer coefficients, $F=(F_1,\dots ,F_n)\in \Z[X_1,\dots ,X_n]^n$. The mysterious Jacobian conjecture (JC)(See [BCW] and [E]), posed first by Keller in 1939 for the integer case and still open even for two dimensional case,  asserts that such a map $F$ is invertible and  has a polynomial inverse with integer coefficients if Jacobian $JF:=\det (\frac {\partial F_i}{\partial X_j})\equiv 1$. K. McKenna and L. van de Dries in [DK] discovered the nice fact that every polynomial surjection of $\Z^n$ is an automorphism of $\Z^n$ that reduces Keller's problem to proving the surjectivity of such maps $F$. 

In this note we present the following.

{\Theorem  Let $F=(P,Q):\C^2\longrightarrow \C^2$ be  a polynomial map with integer coefficients, $F=(P,Q)\in \Z [x,y]^2$ ,  and $JF\equiv 1$. If  the complex plane curve $P=0$ has infinitely many points in $\Z[i]^2$, then $F$  has a polynomial inverse with integer coefficients.}

\medskip
Here, as usual $\Z [i]:=\{ a+bi: a,b\in \Z\}$.

In the proof of Theorem 1 presented on \$ 4 we  will show that   if $F=(P,Q)$ is not invertible, the numbers of points in $\Z[i]^2$ lying on  curves $P=c$, $ c\in\Z[i],$ must be uniformly bounded. This fact will be deduced from Proposition 1, \$ 3, which gives a plane version of van de Dries's result  (see in [CD] and [E]) on the behavior of  integer counterexamples to (JC).

In view of   Siegel's theorem  [Abh. Deutsch. Akad. Wiss. Berlin Kl. Phys.-Mat. 1929, no. 1], which states that there are only finitely many integer points on a curve of genus $g \geq 1$, if the complex plane curve $P=0$ has infinitely many points in $\Z[i]^2$, at least one of it's irreducible components is a rational curve.  As a partial case of the plane Jacobian conjecture it is raised the question {\it whether  a polynomial map $f=(p,q)\in\C[x,y]^2$ with $Jf\equiv c\in\C^*$ is invertible if the curve $p=0$ contains a  rational curve}. However, it was known only that such a map $f$ is invertible if the curve $p=0$ has an irreducible component homeomorphic to $\C$ or if all fibres of $p$ are irreducible and the generic fiber of $p$ is a rational curve (see [R], [LW] and [NN]).

\bigskip

\noindent {\bf 2. Lemma on partial inverse.}
We consider a given $F:\C^n \longrightarrow \C^n$, $(F_1,\dots ,F_n)\in \C[X_1,\dots ,X_n]^n$, with  $F(0)=0$ and $JF\equiv c\in \C^*$. In view of the implicit  function theorem the map $F$ has a unique local analytic inverse $G(Y)=(G_1,\dots , G_n)(Y)$ defined on an open neighborhood $W$ of $0$ for which $G(0)=0$ and 
$$F\circ G(Y)=Y ,\ G\circ F(X)=X\ , \  Y\in W, X\in G(W).\leqno (1) $$
 
Furthermore, the components $G_i$ of $G$ are power series in variables $Y=(Y_1,\dots ,Y_n)$  with complex coefficients and the indenties in (1) hold in the meaning of formal series.

The following lemma, which will be used in the next section, may has an interest in itself.

{\Lemma  Let $1\leq k\leq n $ be fixed and  $L:=\{ y\in \C^n: y_j=0 \ for \ k<j\leq n \}$. If  $G_i(Y_1,\dots ,Y_k, 0, \dots , 0)$, $i=\overline{1,n}$, are polynomials, then

i) the map $g:L\longrightarrow g(L)\subset \C^n$ defined by 
$$g(y_1,\dots ,y_k,0,\dots , 0):= G(y_1,\dots ,y_k,0,\dots , 0)$$is an isomorphism,

ii) the inverse of $g$ is the restriction $f$ of $F$ to the image $g(L)$ and 

iii) If $n=2$ and $k=1$, then  the map $F$ is invertible and $F^{-1}=G$.}

\medskip

\begin{proof}
First, we will prove (i) and (ii). Let $V$ be the connected component of $F^{-1}(L)$ containing the origin $0$. Since $JF\equiv c\in \C^*$,  $F^{-1}(L)$ is a nonsingular algebraic set of dimension $k$, and hence, $V$ is a nonsingular irreducible algebraic set of dimension $k$. Let $F_V$ be the restriction of $F$ on $V$. As assumed, the components $G_i(Y_1,\dots ,Y_k, 0, \dots , 0)$, $i=\overline{1,n}$, are polynomials, the maps $F_V\circ g$ and $g\circ F_V$ are well defined and are regular morphisms. Furthermore,  by (1) we have that $ F_V\circ g=id_L  $ on $L\cap W$ and $g\circ F_V=id_V$ on $V\cap G(W)$.
As $ L$ and $V$ are nonsingular  irreducible algebraic sets of dimension $k$, it follows that
$$ F_V\circ g=id_L ; g\circ F_V=id_V$$
and 
$$g(L)=V ; \  f=F_V.$$
Hence, we get  the conclusions (i) and (ii).

Now, we consider the case $n=2$ and $k=1$. By (i) and (ii) the maps $g: L\longrightarrow \C^2$ and $f: g(L)\longrightarrow \C^2$ are embeddings. Applying the Abhyankar-Moh-Suzuki Theorem on embeddings of the line into the plane (see [AM], [S]), we can find  a new affine coordinate of $\C^2$ in which  $g(L)$ is  a line. Thus, $F$ maps the line $g(L)$ one-to-one to the line $L$, as $f$ is the restriction of $F$ to $g(L)$. Then, $F$ is invertible and $F^{-1}=G$ by the well known fact  (see in [G], [E]) that a non-zero constant Jacobian polynomial map of $\C^2$  is invertible if it sends a line one-to-one into the plane. (This results from an application of Abhyankar-Moh-Suzuki Theorem and the similarity of Newton polygons of components of Jacobian pairs.)
\end{proof}

\bigskip
\noindent {\bf 3. A plane  version of van de Dries's result.}
Van de Dries ( see in [CD] and [E]) observed that if $F\in \Z [X]^2$ is a counterexample  to (JC) and $JF\equiv 1$, then $F$ will maps the lattice $\Z[i]^n$ into a narrow neighborhood of it's exceptional value set $A_F$, namely $$Dist(F(p),A_F)\leq 1 \mbox{ for all } p\in  \Z[i]^n.\leqno (2)$$ Here, $Dist (p,V):=\inf_{q\in V} \max_{i=\overline{1,n}}\vert p_i-q_i\vert$  and the {\it exceptional value set } $A_F$ of  $F$ is the smallest subset $A_F\subset \C^n$ such that the restriction $F:\C^n\setminus F^{-1}( A_F)\longrightarrow \C^n\setminus A_F$ gives a unbranched covering.

The following gives  a version of the mentioned result for the plane case. For points $q=(a,b)\in \C^2$ and a subset $V\subset \C^2$ we define 
$$\hat d (q,V):=\max (\min  \{ \vert a-u\vert: (u,b)\in V \}, \min  \{ \vert b-v\vert: (a,v)\in V \}).$$
Here, as convention $\max \{\emptyset \}:=+\infty .$ Clearly,  $ Dist(q,V)\leq \hat d (q,V)$. 

\medskip
{\Proposition  Let $F=(P,Q)\in \Z [x,y]^2$ with Jacobian $JF\equiv 1$. If $F$ is not invertible, then $$\hat d(F(p),A_F) \leq1 \mbox{ for all } p\in \Z[i]^2\leqno (3) $$}

\begin{proof}
Let $F=(P,Q)\in \Z[x,y]^2$ with $JF\equiv 1$ and suppose that $F$ is not invertible.
Assume the contrary that there is a point $(a,b)\in \Z[i]^2$ such that
$\hat d(F(a,b),A_F)>1$, for instance,
$$\min  \{ \vert P(a,b)-u\vert: (u,Q(a,b))\in A_F \}>1.\leqno (4)$$
Instead of $F$, we consider the map $\bar F=(\bar P,\bar Q)$ given by $\bar F(x,y):=F(x+a,y+b)-F(a,b)$. As $(a,b)\in \Z[i]^2$ and $F\in \Z[x,y]^2$ with $JF\equiv 1$, one can verify that $\bar F\in \Z[i][x,y]^2$, $J\bar F\equiv 1$, $F(0,0)=(0,0)$ and,  like as $F$, $\bar F$ is not  invertible. Furthermore, $A_{\bar F}=A_F-F(a,b)$ and 
$\min  \{ \vert P(a,b)-u\vert: (u,Q(a,b))\in A_F \}=\min  \{ \vert u\vert: (u,0)\in A_{\bar F} \}$.  Hence,  by (4) $$\min  \{ \vert u\vert: (u,0)\in A_{\bar F} \}>1.$$
It follows that there are numbers $r>1$ and $s>0$ small enough such that the set $A_{\bar F} $ does not intersect the box
$B:=\{(u,v)\in \C^2 : \vert u\vert <r; \vert v\vert <s\}$.

Now, let $G(u,v)=(R(u,v),S(u,v))$ be the local inverse of $\bar F$ at $(0,0)$, $G(0,0)=(0,0)$, $F\circ G (u,v)=(u,v)$ on a neighborhood $W$ of $(0,0)$ and $R(u,v)$ and $S(u,v)$ are convergent power series in $W$. Since $\bar F: \C^2\setminus F^{-1}(A_{\bar F}\longrightarrow \C^2\setminus A_{\bar F}$ is a unbranched covering, $B$  is open simple connected set and $B\cap A_{\bar F}=\emptyset$, the local inverse $G$ of $\bar F$ can be extended  analytically over the box $B$. It follows that the power series $R(u,0)$ and $S(u,0)$ are convergent for $\vert u \vert <r$. As $R(u,0)$ and $S(u,0)$ are power series with coefficients in $Z[i]$ and $r>1$,  $R(u,0)$ and $S(u,0)$ must be polynomials in $u$. Hence, by applying Lemma 1 we get  that $\bar F$ is invertible which contradicts to the assumption.
\end{proof}

\noindent{\bf Remark 1.} As shown in [E, p. 262], there exists dominant mappings $F$ of $\C^n$which satisfy van de Dries's estimation. One of such maps is the map $F=(P,Q)$, $P(x,y)=x^6y^4+x^2y$, $Q(x,y)=x^9y^6+3x^5y^3+3x$,  given by Makar-Limanov. For this map  $A_F=\{u=0\}\cup\{u^3-v^2=0\}$, $ P(x,y)=s^2-(xy)^{-2}$ and $Q(x,y)=s^3-(xy)^{-3}$, where $s:=x^3y^2+(xy)^{-1}$. Then,  $ Dist(F(x,y),A_F) =0 $  for $x=0$ or $y=0$ and $ Dist(F(x,y),A_F)\leq Dist(F(x,y),(s^2,s^3)) \leq 1$ for $(xy)^{-1} \leq 1$. Hence, $ Dist(F(x,y),A_F)\leq 1$ for all $(x,y)\in \Z[i]^2$. However, this map does not satisfy the estimation (3). Indeed, for $(x,y)=(1,1)$, $F(1,1)=(3,7)$ and the line $u=3$ intersects $A_F$ at $(3,3\sqrt{3})$ and $(3,-3\sqrt{3})$. So, $\hat d(F(1,1),A_F) >1$, since $\min  \{ \vert 7-v\vert: (3,v)\in A_F \}= 7- 3\sqrt{3}>1$.

\bigskip

\noindent{\bf 4. Proof of Theorem 1.} 
Let $k\in \Z[i]$ and  denote by $I(P,k)$ the set of all points in $\Z[i]^2$ lying on the curve $P=k$, $I(P,k):=\{P=k\}\cap \Z[i]^2.$ Theorem 1 follows directly from the following lemma.

{\Lemma Let $F=(P,Q)\in \Z [x,y]^2$ with Jacobian $JF\equiv 1$. If $F$ is not invertible, then the numbers $\#I(P,k)$, $k\in \Z[i]$, must be uniformly bounded,
$$\max_{k\in \Z[i]} \# I(P,k) <+\infty.$$}

\begin{proof} Assume that $F$ is not invertible. Then, $A_F \neq \emptyset $ and is an algebraic curve in $\C^2$( See for example [J]). Let $Deg A_F$ be the degree of the curve $A_F$ and $deg_{geo.}F$ be the topological degree of $F$, $\deg_{geo.}F=\# F^{-1}(p)$ for generic points $p\in \C^2$.  As $JF\equiv 1$,  $\deg_{geo.}F\geq \# F^{-1}(p)$ for all $p\in \C^2$.

Now, fix a number $k\in \Z[i]$. Let us denote by $L$ the line $\{(u,v)\in \C^2: u=k\}$  and by $D$ the disk $\{ (0,c)\in C^2:  \vert c\vert \leq 1\}$. Then, in view of Proposition 1 we have 
$$I(P,k)\subset F^{-1}(\bigcup_{p\in L\cap A_F}((p+D)\cap \Z [i]^2 ))).$$
It is easy to see that  each of the sets $(p+D)\cap \Z [i]^2$ cannot have more than $5$ points. So, we get
$$\# I(P,k) \leq 5 \deg_{geo.} F \# (L\cap A_F).$$
The number $\# (L\cap A_F)$ cannot be larger than the intersection number of the line $L$ and the curve $A_F$, except for the situation $L\subset A_F$. However, such situation is impossible, since the irreducible components of $A_F$ cannot be isomorphic to a line ( see [C1]).  Hence, we obtain
$$ \# I(P,k) \leq 5 \deg_{geo}F \deg A_F.\leqno(4) $$ 

\end{proof}

\noindent{\bf Remark 2.} 
Examing more carefully the intersection $(L\cap A_F)$ and using results of [C1], one can improve the bound in (4). In fact, we can get  $$ \# I(P,k) < 5 \deg_{geo}F \deg P\frac{\gcd(\deg P,\deg Q)-1}{\gcd(\deg P,\deg Q)}.\leqno(5) $$

\noindent{\bf Remark 3.} One may ask for which curves $h=0$ the set $I(h\circ F,0)$ is finite. If $h=0$ is not such a curve, the curve $h\circ F=0$ must has a branch at infinity $\gamma$  which contains infinitely many integer points. Then, the branch $F(\gamma)$, which is a branch at infinity of $h=0$,   must be very close to  branches at infinity of the curve $A_F$. Using the description of $A_F$ in [C2], we can see that if $(P,Q)\in \Z[x,y]^2$ is a counterexample to (JC), then $\#I(h\circ F,0)<+\infty$ for almost all curves $h=0$ homeomorphic to $\C$. 

\bigskip

To conclude the paper we want to note that both of Theorem 1 and Proposition 1 remain true  if one replaces the ring $\Z[i]$ by  subrings $R$ of the ring of integers of a quadratic number field of the form $Q(i\sqrt {m})$. In fact, in our arguments to prove Theorem 1 and Proposition 1 the only fact of $\Z[i]$  used is that every power series  with coefficients in $\Z[i]$ and convergent radius  larger than $1$ must be polynomial. This fact is also true for power series with coefficients in such a ring $R$ (see [E]).

\bigskip
\noindent{\bf Acknowledgments.} The author wishes to thank Prof.  A. van den Essen   for many valuable suggestions  and useful discussions.

\bigskip

\noindent {\bf References} 

 \noindent [AM]  S.S.  Abhyankar and T.T.   Moh,     {\it  Embeddings of the line in the 
plane},     J.   Reine Angew.   Math.   276 (1975),     148-166.   

\noindent [A]  S. S. Abhyankar,    {\it Expansion Techniques in Algebraic Geometry},    
Tata 
institute of fundamental research,    Tata Institute,    1977.  

\noindent [BCW] H.   Bass,     E.   Connell and D.   Wright,    {\it  The Jacobian conjecture: 
reduction of degree and formal expansion of the inverse},     Bull.   Amer.   Math.   
Soc.   (N.S.) 7 (1982),     287-330.   

\noindent [CD] E. Connell and L. van de Dries, {\it  Injective polynomial maps and the Jacobian conjecture},
J. of Pure and Applied Algebra, 28 (1983), 235-239.

\noindent [C1] Nguyen Van Chau,  {\it 
Non-zero constant Jacobian polynomial maps of $ C^2. $}
J.  Ann.  Pol.  Math.  71,  No. 3 (1999),  287-310.  

\noindent [C2] Nguyen Van Chau, {\it Note on  the Jacobian condition and Non-proper value set}, J.  Ann. Polon. Math. 84, No. 3 (2004), 203-210.

\noindent [DK] L. van de Dries and K. Mckena, {\it Surjective polynomial maps and a remark on the Jacobian problem}, Manuscripta math. , 67(1990), 1-15.

\noindent [E] van den Essen,  Arno,  {\it
Polynomial automorphisms and the Jacobian conjecture}.  (English.  English summary) Progress in Mathematics,  190.  Birkhäuser Verlag,  Basel,  2000. 

\noindent [G] J.  Gwozdziewicz,  {\it Injectivity on one line}.  Bull.  Soc.  Sci.  Lodz 7 (1993),  59-60,  Series: Recherches sur les deformationes XV. 

\noindent [J] Z.  Jelonek, {\it The set of points at which a polynomial map is not proper.} Ann.  Pol.  Math,  58 (1993),  259-266. 

\noindent [LW] Le Dung Trang and C. Weber, {\it  Polynômes à fibres rationnelles et conjecture jacobienne à 2 variables. } C. R. Acad. Sci. Paris Sér. I Math. 320 (1995), no. 5, 581--584. 

\noindent [NN]  W.D. Neumann and P. Norbury, {\it Nontrivial rational polynomials in two variables have reducible fibres.} Bull. Austral. Math. Soc. 58 (1998), no. 3, 501--503.

\noindent [R] M. Razar, {\it  Polynomial maps with constant Jacobian. }
Israel J. Math. 32 (1979), no. 2-3, 97--106.

\noindent [S] M.   Suzuki,     {\it Proprietes topologiques des polynomes de deux variables compleces at automorphismes algebriques  de lespace $\C^2$,    } J.   Math.   Soc.   Japan 26,     2 (1974),     241-257.  
\end{document}